\documentclass[12pt]{article}
\usepackage{amssymb,epsfig,amsfonts,epic,graphics,pictex} 
\begin{document}

\newtheorem{lem}{Lemma}[section]
\newtheorem{prop}{Proposition}
\newtheorem{con}{Construction}[section]
\newtheorem{defi}{Definition}[section]
\newtheorem{coro}{Corollary}[section]
\newtheorem{rem}{Remark}[section]
\newcommand{\hf}{\hat{f}}
\newtheorem{fact}{Fact}[section]
\newtheorem{theo}{Theorem}
\newcommand{\Br}{\Poin}
\newcommand{\Cr}{{\bf Cr}}
\newcommand{\dist}{{\rm dist}}
\newcommand{\diam}{\mbox{diam}\, }
\newcommand{\mod}{{\rm mod}\,}
\newcommand{\compose}{\circ}
\newcommand{\dbar}{\bar{\partial}}
\newcommand{\Def}[1]{{{\em #1}}}
\newcommand{\dx}[1]{\frac{\partial #1}{\partial x}}
\newcommand{\dy}[1]{\frac{\partial #1}{\partial y}}
\newcommand{\Res}[2]{{#1}\raisebox{-.4ex}{$\left|\,_{#2}\right.$}}
\newcommand{\sgn}{{\rm sgn}}

\newcommand{\CC}{\mathbb{C}}
\newcommand{\D}{{\bf D}}
\newcommand{\Dm}{{\bf D_-}}
\newcommand{\RR}{\mathbb{R}}
\newcommand{\NN}{\mathbb{N}}
\newcommand{\HH}{\mathbb{H}}
\newcommand{\Z}{{\bf Z}}
\newcommand{\tr}{\mbox{Tr}\,}
\newcommand{\R}{{\bf R}}
\newcommand{\C}{{\bf C}}

\newenvironment{nproof}[1]{\trivlist\item[\hskip \labelsep{\bf Proof{#1}.}]}
{\begin{flushright} $\square$\end{flushright}\endtrivlist}
\newenvironment{proof}{\begin{nproof}{}}{\end{nproof}}

\newenvironment{block}[1]{\trivlist\item[\hskip \labelsep{{#1}.}]}{\endtrivlist}
\newenvironment{definition}{\begin{block}{\bf Definition}}{\end{block}}

\newtheorem{conjec}{Conjecture}

\newtheorem{com}{Comment}
\font\mathfonta=msam10 at 11pt
\font\mathfontb=msbm10 at 11pt
\def\Bbb#1{\mbox{\mathfontb #1}}
\def\lesssim{\mbox{\mathfonta.}}
\def\suppset{\mbox{\mathfonta{c}}}
\def\subbset{\mbox{\mathfonta{b}}}
\def\grtsim{\mbox{\mathfonta\&}}
\def\gtrsim{\mbox{\mathfonta\&}}

\newcommand{\Poin}{{\bf Poin}}
\newcommand{\Bo}{\Box^{n}_{i}}
\newcommand{\Di}{{\cal D}}
\newcommand{\gd}{{\underline \gamma}}
\newcommand{\gu}{{\underline g }}
\newcommand{\ce}{\mbox{III}}
\newcommand{\be}{\mbox{II}}
\newcommand{\F}{\cal{F}}
\newcommand{\Ci}{\bf{C}}
\newcommand{\ai}{\mbox{I}}
\newcommand{\dupap}{\partial^{+}}
\newcommand{\dm}{\partial^{-}}
\newenvironment{note}{\begin{sc}{\bf Note}}{\end{sc}}
\newenvironment{notes}{\begin{sc}{\bf Notes}\ \par\begin{enumerate}}%
{\end{enumerate}\end{sc}}
\newenvironment{sol}
{{\bf Solution:}\newline}{\begin{flushright}
{\bf QED}\end{flushright}}

\title{On Hausdorff dimension
of some Cantor attractors}
\date{}

\author{G. Levin \\
\small{Dept.\ of Math., Hebrew Univ.}\\
\small{Jerusalem 91904, Israel}\\
\small{\tt levin@math.huji.ac.il}\\
\and
F. Przytycki \\
\small{Inst.\ of Math. of PAN}\\ 
\small{Warsaw 00-950, Poland}\\
\small{\tt F.Przytycki@impan.gov.pl}\\
}
\normalsize
\maketitle
\abstract{We study what happens with the dimension
of Feigenbaum-like attractors of smooth unimodal maps
as the order of the critical point grows}
\section{Introduction}

Let $f$ be a smooth unimodal map of an interval.
We assume that $f$ is infinitely-renormalizable with
stationary combinatorics.
Then $f$ has an attractor $C(f)$ both in
metric and topological senses, which is 
a Cantor set and which is the $\omega$-limit set
of the critical point of $f$.
In this note we consider the following question
motivated by~\cite{bruke},~\cite{novstr}, and~\cite{lsw}: what happens with
the Hausdorff dimension of $C(f)$  
as the order $\ell$ of 
the critical point
grows to infinity? 
We show that it must grow to at least $2/3$.
In the orientation reversing case
(which includes the classical Feigenbaum's one)
we also prove that the Hausdorff dimension has a limit
as $\ell$ tends to infinity, 
this limit is less than
$1$, and it is equal to the Hausdorff dimension of 
an attractor of some limit unimodal dynamics defined in~\cite{lsw}.

Denote by $HD(E)$ the Hausdorff dimension
of a set $E$ in ${\bf R}^n$.

It is well-known~\cite{dmvs} (and follows from convergence
of renormalizations), that 
the Hausdorff dimension $HD(C(f))$ of the 
attractor $C(f)$ of $f$ depends actually
only on the stationary combinatorics
$\aleph$ of the map $f$ and the criticality order $\ell$
of its critical point provided that $\ell$ is
an even integer. It allows us to write $D(\aleph, \ell)=HD(C(f))$
for all smooth $f$ with fixed $\aleph$ and $\ell$.

(Note here that 
once the convergence of renormalizations is established
for all real big enough criticalities $\ell$ 
all results and proofs of the paper hold true for such $\ell$.)

We have a priori:
\begin{equation}\label{apr}
0<HD(\aleph, \ell)< 1.
\end{equation}
\begin{com}\label{first}

(1) If $\ell=2$, then the upper bound in~(\ref{apr}) 
can be strengthened~\cite{grako}: there is a number
$\sigma<1$, such that $HD(\aleph, 2)\le \sigma$
for all combinatorics $\aleph$.

(2) Feigenbaum's case $|\aleph|=2$ with the quadratic critical point
($\ell=2$) has been studied intensively,
see ~\cite{vsk},~\cite{ledmis},
particularly in the framework of
Feigenbaum's universality~\cite{feig0},~\cite{feig1}.
Numerically,
$D(\aleph, 2)=0.538...$, see~\cite{g}.

(3) Although $HD(\aleph, \ell)$ is always positive,
it is not difficult to construct
a sequence of stationary combinatorics $\aleph_n$, such that,
for every $\ell$, $HD(\aleph_n, \ell)\to 0$ as $n\to \infty$.  
For instance, $\aleph_n$ can be defined by
the following first $n-1$ itineraries of the critical value:
$n-2$ times "plus" and one time "minus". 
Then bounds (real or complex)
imply that if $f_n(z)=z^\ell+c_n$ is infinitely-renormalizable with
the stationary combinatorics $\aleph_n$, then
$HD(C(f_n))\to 0$ as $n\to \infty$. 
\end{com}

Note that the number $D(\aleph, 2)$ ($|\aleph|=2$)
as well as the numbers $HD(\aleph_n, \ell)$
(with fixed $\ell$ and big $n$) are less than $2/3$.
\begin{theo}\label{theoKJ} 
For every $\aleph$,
\begin{equation}\label{ieqK}
\liminf D(\aleph, \ell) >{2\over 3}
\end{equation}
as $\ell$ tends to infinity along the even integers.
\end{theo}
To state our result about the upper bound, we need
to introduce some notions.

{\it Non-symmetry.} 
For a unimodal map $f$ with a single critical point at $c$,
denote by $I_f$ the involution map
defined in a neighborhood of $c$ by $I_f:x\mapsto \hat x$, where
$I_f(c)=c$, and otherwise $I_f(x)$ is the unique
$\hat x\not=x$, such that $f(x)=f(\hat x)$.
If $f$ is of the form $|E(x)|^\ell$,
where $\ell>1$ and $E$ is a $C^2$-diffeomorphism,
then $I_f$ is also $C^2$, and $I'_f(c)=-1$.
The {\it non-symmetry} $N(f)$ of $f$ is said to be the 
number $N(f)=|I''_f(c)/2|$. It is easy to check that
$N(f)=|E''(c)/E'(c)|$.         

{\it Orientation reversing combinatorics}
of an infinitely-renormalizable unimodal map $f$
is such stationary combinatorics $\aleph$, that
the rescaling factor of the renormalization
is negative. In other words, the maps $f$ and 
$f^{|\aleph|}$ have at the critical point of $f$
different type of extrema (maximum and minimum).
Examples: $|\aleph|=2,3$; more generally, $\aleph_n$ ($n\ge 1$) defined in 
Comment~\ref{first}(3).

For a combinatorial type $\aleph$ and an even integer $\ell$,
denote by $H_{\aleph, \ell}$ the unique universal unimodal map
normalized so that $H_{\aleph, \ell}:[0,1]\to [0,1]$
and $H_{\aleph, \ell}(0)=1$ (see next Section for complete
definition). It is shown in~\cite{lsw}, that
the sequence  $\{H_{\aleph, \ell}\}_{\ell}$ converges uniformly
to a unimodal map $H_\aleph: [0,1]\to [0,1]$.

We prove in Lemma~\ref{dom} that 
{\it if the combinatorial type $\aleph$ reverses orientation, then
the sequence of non-symmetries
$N(H_{\aleph, \ell})$, $\ell=2,4,...$, is uniformly bounded.}

\begin{theo}\label{cont}
For a given combinatorial type $\aleph$, assume that
the sequence of non-symmetries 
$N(H_{\aleph, \ell})$, $\ell=2,4,...$, is uniformly bounded.
Then the Hausdorff dimension of the attractor
is continuous at $\ell=\infty$: there exists
\begin{equation}\label{limhd}
\lim_{\ell\to \infty} D(\aleph, \ell)=HD(C(H_{\aleph}))<1.
\end{equation}
Consequently, ~(\ref{limhd}) holds when $\aleph$ reverses
orientation.
\end{theo}
\begin{com}
It is not clear if the non-symmetry $N(H_{\aleph, \ell})$
is uniformly bounded in $\ell$ for any type $\aleph$.
\end{com}

The proof of Theorems~\ref{theoKJ}-~\ref{cont}
is based on recent results 
of~\cite{lsw}: see next Sect. where we reduce the statements
to Theorem~\ref{cover}. 

(Note however that in the proof of the lower
$2/3$-bound we use only a part of the main result
of~\cite{lsw}, namely, the compactness (Theorem 4 in~\cite{lsw}).)

In turn, to prove Theorem~\ref{cover} we use some results
of~\cite{mu},~\cite{mcmII}, see Sect.~\ref{back}.

From now on, we fix the type $\aleph$. Denote $p=|\aleph|$.

{\bf Acknowledgment.} The first author thanks Benjamin Weiss for a helpful
discussion. The second author acknowledges Polish KBN grant 2PO3A 034 25.

\section{Reduction to fixed-point maps}\label{red}
\subsection{Universal maps}
For every real number $\ell>1$,
we consider a unimodal map
$g_{\ell}: [-1, 1]\to [-1, 1]$ with the
critical 
point at $0$
of order $\ell$. More precisely, $g_\ell$ is assumed to be
in the following form:
$g_{\ell}(x)=E_\ell(|x|^\ell)$, where $E_\ell :\: [0,1] \rightarrow
\RR$ is a $C^2$-diffeomorphism onto its image. 
The map $g=g_\ell$ is normalized so that $g_\ell(0)=1$.
It is further assumed to be infinitely renormalizable
with the fixed combinatorial order type $\aleph$ and to satisfy the fixed
point equation:

\begin{equation}\label{equfeig}
 \alpha g^{|\aleph|}(x) = g(\alpha x) \; .
\end{equation}

with $|\alpha|>1$. 
By renormalization theory, see~\cite{su}, a fixed point $g_{\ell}$ for
any $\ell>1$ can be represented as $E_{\ell}(|x^\ell|)$
with $E_{\ell}$ which is
a diffeomorphism in Epstein class (i.e.
a diffeomorphism $E$ of a real interval $T'$ onto another real
interval $T$ such that
the inverse map $E^{-1}: T\to T'$ extends to a univalent
map $E^{-1} :\: (\CC \setminus \RR) \cup T \rightarrow
(\CC\setminus\RR) \cup T'$).

It will be useful to deal with another unimodal map $H_\ell$,
which is related to $g_\ell$ as follows:
$H_\ell(x)=|g_\ell(x^{1/\ell})|^\ell=|E_\ell(x)|^\ell$, $0\le x\le 1$.
Then $H_\ell$ is a unimodal map of $[0, 1]$ into itself,
with a strict minimum attained at some
$x_\ell\in (0, 1)$. It also satisfies the equation:

\begin{equation}\label{equfeigH}
 \tau H^{|\aleph|}(x) = H(\tau x) \; .
\end{equation}

with $\tau=|\alpha|^\ell$.


We denote by $C(g_\ell)$ and $C(H_\ell)$ the attracting Cantor sets
of the maps $g_\ell: [-1, 1]\to [-1, 1]$ and
$H_\ell: [0, 1]\to [0, 1]$ respectively.
Clearly, $HD(C(g_\ell))=HD(C(H_\ell))$.
Indeed, $E$ conjugates $H=H_\ell$ to $g$ restricted to $[g(1),1]$, 
therefore it maps $C(H)$ to $C(g)$ and is a diffomorphism between  
neighbouhoods of these sets.

Assume now that the order $\ell$ is an even integer.
Then the equation ~(\ref{equfeig}) with the normalization as above
does have a unique solution, for every fixed $\ell$ and $\aleph$,
see~\cite{su},~\cite{mcm}. Consequently, $H_\ell=|g_\ell(x^{1/\ell})|^\ell$
is the unique solution of
~(\ref{equfeigH}) with the normalization as above.

In what follows, $\ell$ is an even integer, and $H_\ell$ denotes
this unique solution of~(\ref{equfeigH}), with its own scaling constant
$\tau_{\ell} > 1$. (Remind that the type $\aleph$ is fixed.)

\subsection{Limit dynamics}
The following result is proved in~\cite{lsw} (even for real $\ell$),
see Theorems 1-2
and Proposition 3 there:
\begin{theo}\label{theolimit}

The sequence of maps $H_{\ell}$ converges as  
$\ell\rightarrow \infty$, uniformly on $[0,1]$, to a unimodal function
$H=H_\infty$, which satisfies the following properties:

\begin{enumerate}    
\item
$\lim_{\ell\rightarrow \infty} \tau_{\ell} = \tau>1$ exists.
and $H, \tau$ satisfy the fixed point equation
$\tau H^{p} (x) = H(\tau x)$
for every $0 \leq x \leq \tau^{-1}$. Here (as always)
$p=|\aleph|$. 
\item
$H$ has analytic continuation to the union of two topological disks 
$U_{-}$ and $U_{+}$ and this analytic continuation will also be
denoted by $H$. 
\item
For some $R>1$, $H$ restricted to either $U_+$ or $U_-$ is a covering
(unbranched) of the punctured disk $V := D(0, R) \setminus \{0\}$ and
$\overline{U_+ \cup U_-} \subset D(0,R)$.  
\item $U_{\pm}$ are both symmetric with respect to the real axis and their
closures intersect exactly at $x_0$; 
$[0,x_0)\subset U_-$, $(x_0,1]\subset U_+$.  
\item Each $H_{\ell}$ extends to complex-analytic map
defined in $U_-\cup U_+$;
this sequence
of analytic extensions converges to $H$,
as $\ell_m\to \infty$, uniformly
on every compact subset of $U_-\cup U_+$.
\item For any two open intervals $I, J$ of the real axis,
if $0\notin J$ and $H:I\to J$ is one-to-one,
then the branch $H^{-1}:J\to I$ extends to a univalent
map to the slit complex plane $({\bf C}\setminus {\bf R})\cup J$
(this follows from the same property for $H_\ell$ with $\ell$ finite)
\item
The mapping $G_\infty(x) := H^{p-1}(\tau^{-1} x)$ fixes $x_0$ 
and $G_\infty^2$ has the following power
series expansion at $x_0$:
\[ G_\infty^2(x) = x - a (x-x_0)^3 + O(|x-x_0|^4) \]
with $a> 0$.  
\item For each $\ell$, the mapping
$G_\ell:=H_\ell^{p-1}(\tau_\ell^{-1} x)$ fixes the critical point
$x_\ell$ of $H_\ell$, $G'_\ell(x_\ell)=\pm 1/\tau_\ell^{1/\ell}$,
and $G_\ell$ converge to $G_\infty$ uniformly in a (complex)
neighborhood of $x_0$.
\item The unimodal map map $H: [0,1]\to [0,1]$ has a unique attractor
$C(H)$, which (as for finite $\ell$) is the closure
of iterates of the critical point. 
\end{enumerate}
\end{theo}
\subsection{The reduction}
Since we know already that $HD(C(f))$ depends merely
on $\aleph$ and $\ell$, Theorems~\ref{theoKJ}-~\ref{cont}
are covered by
the following statement 
\begin{theo}\label{cover}
The following holds.

(a)
\begin{equation}\label{coverK}
\liminf_{\ell\to \infty} HD(C(H_{\ell}))\ge HD(C(H_{\infty})).
\end{equation}

(b) 
\begin{equation}\label{coverinf}
{2\over3}< HD(C(H_{\infty}))<1;
\end{equation}

(c) if the non-symmetries $N(H_{\ell})$
are uniformly bounded as $\ell\to \infty$, then 
the Hausdorff dimension is continuous at infinity:
\begin{equation}\label{covercont}
\lim_{\ell\to \infty} HD(C(H_{\ell}))=HD(C(H_{\infty})).
\end{equation}
\end{theo}
The rest of the paper is devoted to the proof
of this statement.

\section{Background in dynamics}\label{back}
We prove Theorem~\ref{cover}
by reducing it finally to known statements
about infinite conformal iterated function systems (c.i.f.s.)~\cite{mu}
and asymptotics near parabolic maps~\cite{mcmII}, which are given here.
\subsection{C.I.F.S.}
We follow~\cite{mu} restricting ourself to dimension one.
Let $X$ be a closed real interval, and $\sigma$ be
a positive continuous function on $X$, which defines
a new metric $d\rho=\sigma dx$ on $X$.
Let $I$ be a countable index set, $|I|>1$, and
let $S=\{\phi_i: X\to X, i\in I\}$ be a collection
of injective uniform contractions w.r.t. the metric $\rho$:
there is $\lambda<1$, such that 
$\rho(\phi_i(x), \phi_i(y))\le \lambda \rho(x,y)$
for all $i$ and all $x,y$.
For every finite word $w=w_1...w_n$, denote
$\phi_w=\phi_{w_1}\circ...\circ \phi_{w_n}$.
(Note that the metric $\rho$ can be replaced by the
Euclidean one by replacing $\phi_i$ 
by $\phi_w$, where $w$ runs over all finite
words of some fixed length $n$, s.t. $\lambda^n||\sigma||<1$.)
For any infinite word of symbols 
$w=w_1w_2...w_j...$, $w_j\in I$, denote
$w|n=w_1w_2...w_n$. The limit set $L$ of $S$ is 
$L=\cup_{w\in I^\infty} \cap_{n=1}^\infty \phi_{w|n}(X)$.
The system $S$ is said to be conformal if:

(a) $\phi_i(Int(X))\subset Int(X)$ and  
$\phi_i(Int(X))\cap  \phi_j(Int(X))=\emptyset$
for all indexes $i\not=j$.

(b) There is an open set $Y\supset X$, such that
all maps $\phi_i$ extend to $C^{1+\epsilon}$ diffeomorphisms
of $V$ into $V$.

(c) There is $K\ge 1$, such that
$|D\phi_{w}(y)|\le K |D\phi_{w}(x)|$ for every finite
word $w$ and all $x,y\in Y$, where
$D\phi_{w}(x)$ means the derivative w.r.t. the metric $\rho$

The main object of our interest is the
Hausdorff dimension of the limit set $L$.
Note that it is the same w.r.t. the metric $\rho$
as w.r.t. the standard Euclidean metric.

For every integer $n\ge 1$ and every $t\ge 0$ define
$p_n(t)=\sum_{w} ||D\phi_w||^t$
where $w$ runs over all words of length $n$, and 
$||.||$ means the sup-norm.
Consequently,
$P(t)=\lim_{n\to \infty}{1\over n}\log p_n(t)$
is called the pressure of $S$ at $t$.
The parameter $\theta=\theta_S$ of the system
is defined as $\inf\{t: p_1(t)<\infty\}$.
\begin{theo}\label{mu1}

1. (see~\cite{mu}. Prop. 3.3) $P(t)$ is non-increasing on $[0,\infty)$,
strictly decreasing, continuous and convex on $[\theta,\infty)$.

2. (see~\cite{mu}, Thm. 3.15)  $HD(L)=\sup\{HD(L_F):  F\subset I \ \ is \ \ finite\}=
\inf\{t: P(t)\le 0\}$;
if $P(t)=0$ then $t=HD(L)$.

3. If the series
$p_1(\theta)$ diverges, then $P(HD(L))=0$ and
$\theta<HD(L)$.
\end{theo}

(Note that 3 follows directly from 1-2.)

The system with $P(t)=0$ is called {\it regular}. 
The system is regular if and only if there is a $t$-{\it conformal measure},
i.e. a probability measure $m$ such that $m(L)=1$
and for every Borel set $A\subset X$ and every $i\in I$
$m(\phi_i(A))=\int_A |D\phi_i|^t dm$ and 
$m(\phi_i(X)\cap \phi_j(X))=0$  for all $i\not=j$ from $I$.
\subsection{Dominant convergence and forward Poincar\'e series}
Here we follow~\cite{mcmII} adapting the statements sligthly for our
applications.

Let $f_n:U\to {\bf C}$ be a sequence of holomorphic maps
which converges uniformly in a topological disk $U$ 
of the plane to a holomorphic map $f:U\to {\bf C}$.
Assume that $c_n\to c\in U$, 
and the following expansions hold:
$f_n(z)=c_n+\lambda_n(z-c_n)+b_n(z-c_n)^2-a_n(z-c_n)^3+...$,
where $0<\lambda_n<1$, $b_n, a_n\in {\bf R}$, and
$f(z)=z-a(z-c)^3+...$, where $a>0$, i.e.
$f$ is parabolic with two (``real'') attracting petals at $c$.
(In particular, $b_n\to 0$ and $a_n\to a$.)
Then $f_n$ is said to converge to $f$ {\it dominantly}, if
there is $M>0$ such that
$|b_n|\le M|\lambda_n-1|$ for all $n$. 

For every $g=f_n$ and $t>0$ define the (forward) Poincar\'e
series $P_t(g,x)=\sum_{i\ge 0} |(g^i)'(x)|^t$,
and, for any open set $V\subset U$, define 
$P_t(g,V,x)=\sum_{g^i(x)\in V} |(g^i)'(x)|^t$.
We say the Poincare series for  $(f_n, t_n)$ converge uniformly,
if, for any compact set $K$ ($c\notin K$) 
in an attracting petal of $f$, 
and any $\epsilon>0$ there exists a neighborhood $V$ of $c$, such that
$P_{t_n}(f_n,V,x)<\epsilon$ for all $n$ large enough 
and all $x\in K$. We will need
\begin{theo}\label{poin}
Let $f_n, f$ be as above, and $t_n\to t>2/3$.
If $f_n\to f$ dominantly, then the Poincare series
for  $(f_n, t_n)$ converge uniformly.
\end{theo}
This is a particular case of Theorem 10.2 proven in~\cite{mcmII}.
For completeness,
we give a short proof of Theorem~\ref{poin}, see Appendix.

\section{Proof of Theorem~\ref{cover}} 
\subsection{Presentation system for the Cantor attractor}
We repeat (with modifications) 
a construction from~\cite{lsw} (cf.~\cite{ledmis},~\cite{CEp222}),
which is crucial for our proof.
Let $H$ be either one of $H_\ell$ or the limit map 
$H_\infty$. Consequently, let $G$ be either the corresponding
$G_\ell$ or $G_\infty$.
We construct the presentation system for the attractor
$C(H)$, which is an infinite iterated function 
system $\Pi$ on an interval $I$ 
so that $C(H)\cap I$ is (up to a countable set)  
the limit set of $\Pi$. 
Moreover, this picture converges, as $\ell\to \infty$,
to the corresponding picture of the limit map. 

Denote $c_j=H^{j-1}(0)$, $j\ge 0$, 
the $j$-iterate of the critical point $c_0$ of $H$
(i.e., $c_0=x_\ell$ for $H=H_\ell$ and $c_0=x_0$ for
$H=H_\infty$).
Let $I=[c_p,c_{2p}]$. Then we define a sequence of maps 
$\psi_{k,m}:I\to I$, $k=1,2,...$, $m=1,2,...,p-1$,  as follows.
Let $H^{-(p-m)}: [c_p, c_{2p}]\to [c_m,c_{p+m}]$ 
denote corresponding one-to-one branch of $H^{-(p-m)}$.
Then set
\begin{equation}
\psi_{k,m}=G^k\circ H^{-(p-m)}.
\end{equation}
\begin{lem}\label{ref}

(a) 
$$I_{k,m}:=\psi_{k,m}(I)=[c_{p^k m}, c_{p^k(p+m)}]\subset I.$$
The intervals $I_{k,m}$ are pairwise disjoint.

(b) Let $L$ be the limit set of the system $\{\psi_{k,m}\}$
(in other words, $L$ is the set of non-escaping points
of the inverse maps $\psi_{k,m}^{-1}:I_{k,m}\to I$).
Then the closure $\overline L=L\cup P$,
where $P$ is a subset of pre-images of the critical point
$c_0$, and
$$\overline L=C(H)\cap I.$$
\end{lem}
\begin{proof} From the functional equation for $H$,
$G(c_j)=c_{pj}$, $j\in \Z$, where $c_j$, for $j<0$
is an $H^j$-preimage of $c_0$. The rest follows.
\end{proof}
Denote by $\Pi_\ell=(\psi^{(\ell)}_{k,m})_{k,m}$, resp. 
$\Pi_\infty=(\psi^{(\infty)}_{k,m})_{k,m}$, the
presentation system of $H_\ell$, resp. $H_\infty$.

The notation $B(E)$ stands for the round disk which
is based on an interval $E\subset {\bf R}$ 
as a diameter.
\begin{lem}
Let $\Pi=\{\psi_{k,m}:I\to I_{k,m}\}_{k,m}$ 
be either $\Pi_\ell$ or $\Pi_\infty$.

(1) There exists a fixed open interval $J$, which contains
$I$ for all $\ell$ largh enough (including $\ell=\infty$),
such that each $\psi_{k,m}$ extends to a univalent
map $\psi_{k,m}:B(J)\to B(J_{k,m})$, where 
$J_{k,m}=\psi_{k,m}(J)$ are pairwise disjoint
intervals properly contained in $J$.

Therefore, there is $\lambda<1$ (dependent only on the
type $\aleph$), such that
$||D\psi_{k,m}||_\rho<\lambda$, for all $k,m$,
and $\ell\le \infty$ large enough, where
$||D\psi_{k,m}||_\rho$ denotes the supremum
on the interval $I$
of the derivative of $\psi_{k,m}$ in the hyperbolic
metric $\rho$ of $B(J)$.

(2) $\Pi$ (with the metric $\rho$
restricted to the closed subinterval
$I$ of $J$)
is an infinite conformal iterated function system, such that:

(a) $\theta_{\Pi_\ell}=0$ for $\ell<\infty$;

(b) $\theta_{\Pi_\infty}=2/3$,  $P(\theta_{\Pi_\infty})=\infty$;

(c) $\Pi_\ell$, $\ell\le \infty$, is regular.
\end{lem}
\begin{proof}
(1) follows from Theorem~\ref{theolimit} (7), and from 
another representation of the maps of the system:
$\psi_{k,m}=H^{-1}\circ \tau^{-k}\circ H^{-(p-m-1)}$
which is a consequence of the eq. $H\circ G=\tau^{-1}\circ H$. 
(2a) is immediate because $c_0$ is the attracting fixed point
of $G$ for finite $\ell$.

(2b)-(2c):
since $G=G_\infty$ has a neutral fixed point
with two attracting petals, and $\psi_{k,m}'(x)=
(G^k)'(H^{-(p-m)}(x))(H^{-(p-m)})'(x)$, 
we obtain the following asymptotics, as $k\to \infty$, 
for the 
presentation
system: $|\psi_{k,m}'(x)|/k^{-3/2}\to a_m(x)$
where, for fixed $m=1,...,p-1$, the function
$a_m(x)$ is continuous and positive on $I$.
It follows from here that the critical exponent
$\theta$ of the system is $\theta=2/3$.
Thus, $p_1(\theta)=\infty$ for all $\ell\le \infty$. Hence,
by Theorem~\ref{mu1}, 
the system $\{\psi_{k,m}\}$ is regular.

\end{proof}

\subsection{Hausdorff dimension for the limit map}
As a corollary, we obtain Theorem~\ref{cover},
(a)-(b):
\begin{coro}\label{ab}

(1) $2/3<HD(C(H_\infty))<1$, 

(2) $\liminf_{\ell\to \infty} HD(C(H_\ell))\ge HD(C(H_\infty))>{2\over 
3}$.
\end{coro}
\begin{proof}
Denote $H=H_\infty$.
Since $H$ is regular and $P(2/3)=\infty$, then
$HD(C(H))>2/3$.
On the other hand, the Lebesgue measure of
$I\setminus \cup_{k,m} I_{k,m}$ is positive. 
Therefore (~\cite{mu}, Theorem 4.5),  $HD(C(H))=HD(C(H)\cap I)=HD(L)<1$.

(2) follows from Theorem~\ref{mu1}:
for every $\delta>0$, there is a finite subsystem $F_\infty$
of $\Pi_\infty$ with the Hausdorff dimension
of its limit set at least $HD(C(H_\infty))-\delta$.
Since corresponding finite subsystem $F_\ell$ converges
to $F_\infty$ as $\ell\to \infty$, then
the Hausdorff dimension of the limit set of $F_\ell$
is at least   $HD(C(H_\infty))-2\delta$, for all $\ell$
large enough. The result follows.
\end{proof}

\subsection{Non-symmetry and dominant convergence}
It remains to prove Theorem~\ref{cover} (c).

Denote $\epsilon=1$ or $2$ depending on
whether $G_\infty'(x_0)=1$ or $-1$.
\begin{lem}\label{dom}
1. The sequence $G^\epsilon_\ell$ converges to
$G_\infty^\epsilon$ dominantly if and only if
the sequence of non-symmetries $N(H_\ell)$
is bounded.

2. If the combinatorics reverses orientation, then
$G_\ell^2$ converges dominantly to $G_\infty^2$,
and the non-symmetries $N(H_\ell)$ are uniformly bounded.
\end{lem}
\begin{proof}
Let $H=H_\ell$ and $G=G_\ell$, $\tau=\tau_\ell$, and $I=I_H$.
We have: $H(G(I(x)))=\tau^{-1} H(I(x))=\tau^{-1} H(x)=H(G(x))$,
i.e. $I\circ G=G\circ I$. The latter equation gives us:
$|(G^\epsilon)''(x_\ell)|=N(H)\lambda (1-\lambda)$,
where $\lambda=\lambda_\ell=(G^\epsilon)'(x_\ell)\in (0,1)$.
This implies 1.

To prove 2, notice that the combinatorics reverses
orientation if and only if $G_\infty'(x_0)=-1$.
Then we get the dominant convergence,
because $|(G^2)''(x_\ell)|=|G''(x_\ell)| |\lambda| (1-|\lambda|)$
and $G''(x_\ell)=G_\ell''(x_\ell)$ converges
to the number $G_\infty''(x_0)$, as $\ell\to \infty$.
(One can also refer
formally to~\cite{mcmII}, Proposition 7.3.)
\end{proof}

\subsection{Conformal measures of the presentation systems}
Remind that 
$\Pi_\ell=(\psi^{(\ell)}_{k,m}:I^\ell\to I_{k,m}^\ell)_{k,m}$, resp. 
$\Pi_\infty=(\psi^{(\infty)}_{k,m}:I^\infty\to I_{k,m}^\infty)_{k,m}$, the
presentation system of $H_\ell$, resp. $H_\infty$.
We know that $\Pi_\ell, \Pi_\infty$ are regular. 
Denote by $\mu_\ell$, resp. $\mu_\infty$,
the unique probability $h_\ell$-conformal, 
resp. $h_\infty$-conformal, measure 
of $\Pi_\ell$, resp. $\Pi_\infty$,
where $h_\ell=HD(C(H_\ell)\cap I^\ell)=HD(C(H_\ell))$,
$h_\infty=HD(C(H_\infty)\cap I^\infty)=HD(C(H_\infty))$.
(Notice that the measures have nothing to do with conformal
measures of $H_\ell, H_\infty$, because the dynamics
are completely different.)
Since any regular system has a unique conformal measure,
to prove that $h_\ell\to h_\infty$,
it is enough to prove that a weak limit $\nu$ of a subsequence 
of $\mu_\ell$ is a conformal measure of $\Pi_\infty$.
For this to be true, it is enough to check that
the support of $\nu$ is contained in the limit set
$L_\infty$ of $\Pi_\infty$. Note that by Lemma~\ref{ref}(b),
the set $\overline L_\infty\setminus L_\infty$ is countable.
Therefore, it is enough to prove that $\nu$ has no atoms.
Thus Theorem~\ref{cover}(c) follows from
\begin{lem}
If the non-symmetries $N(H_\ell)$
are uniformly bounded, then the measure $\nu$ has no atoms.
\end{lem}
\begin{proof} Let the point $a\in supp(\nu)=\overline L_\infty$,
where  $L_\infty$ is the limit set of $\Pi_\infty$,
be an atom of $\nu$. Then there is $\sigma>0$
such that for all $r>0$ small enough
$\mu_\ell(B(a,r))>\sigma$ along a subsequence
of $\ell$'s.
Since $\psi_{k,m}$ are
uniform contractioncs and the measures are probabilities, one sees that
$a\in \overline L_\infty \setminus L_\infty$, i.e., afterall,
one can assume that $a=x_0$.
Now 
$\mu_\ell(B(x_0, r))\le 
\sum_{I^\ell_{k,m}\cap B(x_0,r)\not=\emptyset}\int_{I^\ell}
|D\psi_{k,m}^{(\ell)}|^{h_\ell}d\mu_\ell\le
C\sum |(G_\ell^k)'(y_{\ell,m})|^{h_\ell}$,
for some fixed $C>0$, some points $y_{\ell, m}$
from  a fixed compact set $K$, $x_0\notin K$
(if $\ell$ is big enough), and
the latter sum runs over such $k$ that
$G_\ell^k(y_{\ell,m})\in B(x_0, r')$, where
$r'\to 0$ as $r\to 0$. 
Then a contradiction
follows directly from Lemma~\ref{dom} and Theorem~\ref{poin}
(note that $t>2/3$ by Corollary~\ref{ab}(2)).
\end{proof}

\section{Appendix: proof of Theorem~\ref{poin}}
{\bf 1.} If $h_n\to h$ is a sequence of injective holomorphic maps
in a fixed neighborhood of $c$, which converges to an injective $h$
uniformly, then the Poincar\'e series
for  $(f_n, t_n)$ converge uniformly iff
the Poincar\'e series
for  $(h_n\circ f_n\circ h_n^{-1}, t_n)$ converge uniformly.
In particular, one can assume that $c_n=c=0$.

{\bf 2.} (see Theorem 7.2 of~\cite{mcmII}).
Let $h_n(z)=z-B_nz^2$, where $B_n=b_n/(\lambda_n(\lambda_n-1))$.
Since $|b_n|\le M|\lambda_n-1|$ for all $n$,there is a subsequence 
of $h_n$ as in Step {\bf 1}. On the other hand,
$h_n\circ f_n\circ h_n^{-1}(z)=\lambda_n z + O(z^3)$.
It means one can assume that
$f_n(z)=\lambda_n-a_nz^3+...$ where $a_n\to a>0$,
$0<\lambda_n<1$ and $\lambda_n\to 1$.


{\bf 3.} For $f_n$, make a change 
$z=\hat h_n(w)=d_n w^{-1/2}$, where $w\in F=\{w: Re(w)>R_0\}$ and
$d_n= (\lambda_n^3/(2a_n))^{1/2}$.
For $g_n=\hat h_n^{-1}\circ f_n\circ \hat h_n$, it holds
$g_n(w)=\sigma_n w + 1 + \alpha_n(w)$, where
$\sigma_n=\lambda_n^{-2}>1$ and $\sigma_n\to 1$,
$\alpha_n$ converge
uniformly in $F$ to the corresponding 
$\alpha$ for $g=\hat h^{-1}\circ f\circ \hat h$,
$\hat h=\lim \hat h_n$, and  $\alpha_n(w)=O(|w|^{-1/2})$,
$\alpha(w)=O(|w|^{-1/2})$. 

To deal with $g_n^i$, we prove the following 
simple Claim. This is weaker than Theorems 8.1-8.3 of~\cite{mcmII},
but still enough for our needs.

{\bf Claim 1}:
{\it For every $\delta>0$ there is $R_\delta>R_0$ 
and, for every $n$, there is $1+\delta$-quasiconformal map
$\phi_n$ of the plane that fixes $0, 1$, and $\infty$, such that
$\phi_n^{-1}\circ g_n\circ \phi_n=T_n$, where
$T_n(w)=\sigma_n w+1$, for
$Re(w)>R_\delta$. Passing to a subsequence, one can assume that
$\phi_n\to \phi$, so that $\phi^{-1}\circ g\circ \phi=T$,
$T(w)=w+1$.}

{\bf Proof.} 
Fix $\delta>0$. Denote $\Pi(R_1, R_2)=\{w: R_1< Re(w)< R_2\}$.
Then $|\alpha_n(w)|$ and 
$|\alpha_n'(w)|\le \sup\{|\alpha_n(t)|: |t-w|<1\}$
are uniformly arbitrary small as $w\in L:=\{\Re(w)=R_\delta\}$
and $R_\delta\to \infty$.
Therefore, all $\sigma_nw$ can be joined to
$z(w):=\sigma_nw+1+\alpha_n(w)$ by disjoint intervals 
$I(w)$ in the strip between 
$\sigma_n L$ and $z(L)$. The mapping $\phi_n$,
which is affine on
each interval $[\sigma_nw, \sigma_nw+1]$ onto $I(w)$
together with the identity on $\Pi(R_\delta, \sigma_n R_\delta)$,
is $1+\delta$ quasi-conformal on $\Pi(R_\delta, \sigma_n R_\delta+1)$.
Then we extend $\phi_n$ to $Re(w)> \sigma_n R_\delta+1$
by the (conformal) dynamics of $g_n$, $T_n$, and define it identity
on the rest of the plane.

{\bf Claim 2.} {\it For every real $p>1$,
there is $M$ such that
${|(T_n^i)'(w)|\over |T_n^i(w)|^p}\le M i^{-p}$
for all $i,n$, and all $w>1$.}

Indeed, denote $C(i,n)=\sigma_n^i$. Consider any subsequence
$(i_j,n_j), j\to \infty$.
If $C(i,n)$ is bounded from above along this subsequence,
then applying as in~\cite{mcmI}, Sect.6, the inequality
between arithmetic and geometric means, we can write
$T_n^i(w)=\sigma_n^i w + (1+\sigma_n+...+\sigma_n^{i-1})\ge
(i+1)w^{1/(i+1)}\sigma_n^{i/2}\ge  C(i,n)^{1/2} i$, so that
${|(T_n^i)'(w)|\over |T_n^i(w)|^p}\le C(i,n)^{1-p/2} i^{-p}=
O(i^{-p})$ along the subsequence. If now
$C(i,n)\to \infty$  along $(i_j,n_j)$
(and $\sigma_{n_j}\to 1$), then
${|(T_n^i)'(w)|\over |T_n^i(w)|^{p}}
={|\sigma_n^i|\over |\sigma_n^iw + (\sigma_n^i-1)/(\sigma_n-1)|^p}\sim
C(i,n)|\sigma_n-1|^p/C(i,n)^p
\sim 
(\log C(i,n))^{p}/C(i,n)^{p-1} i^{-p}=o(i^{-p})$. 

{\bf 4.} From Steps 1-2, Claim 1, and  Koebe distortion theorem,
it follows that
it is enough to prove the theorem assuming that
the compact $K$ is a point $x$, which moreover
lies on an attracting direction of $f$,
and small neighborhood $V$ can be replaced by big indexes.
We have:
$|(f_n^i)'(x)|=K|(g_n^i)'(w)|/|g_n^i(w)|^{3/2}$, where
$K>0$ and $w>R$ depend only on $x>0$.
Thus we need to show that, if $t_n\to t>2/3$,
for a given $w>0$ close enough to $+\infty$,
for any $\epsilon>0$ there exists an index $i_0$, such that
$S(g_n,i_0,t_n):=\sum_{i\ge i_0} 
|(g_n^i)'(w)/g_n^i(w)^{3/2}|^{t_n}<\epsilon$
for all $n$ large enough. 
Claim 2 (with $p=3/2$) 
implies immediately that this is true for $g_n=T_n$.

To handle $S(g_n,i_0,t_n)$ in general, we 
compare it with $S(T_n,i_0,t_n)$ 
and proceed similar to~\cite{mcmII}, Sect.10.
Due to Koebe distortion theorem, one can replace
the derivative by the ratio of diameters.
By Claim 1, the change of the diameters when passing
from $g_n$ to $T_n$ is H\"older with 
the exponent arbitrary close to $1$.
Then we apply Claim 2 with $p$ arbitrary close to $3/2$.


\end{document}